\pdfoutput=1
\documentclass{article}
\usepackage[margin=25mm]{geometry}
\usepackage{amsmath}
\usepackage{amsfonts}
\usepackage{amssymb}
\usepackage{graphicx}
\usepackage[utf8]{inputenc}
\usepackage[T2A]{fontenc}
\usepackage[ukrainian,russian,english]{babel}
\usepackage[unicode, pdftex]{hyperref}
\usepackage{xcolor}
\usepackage[all]{xy}

\newtheorem{theorem}{Theorem}

\title{
Structures of optimal discrete gradient vector fields on surface with one or two critical cells  }

\author{Svitlana Bilun, Maria Hrechko, Olena Myshnova and Alexandr Prishlyak}
\begin{document}

\maketitle
\begin{abstract}
   We describe all possible  structures of discrete vector field  (discrete Morse functions) with minimal number of critical cells on the regular CW-complex for  the 2-disk (1 cell), the 2-sphere (2 cells), the cylinder (2 cells) and Mobius band (2 cells). 
\end{abstract}
\section*{Introduction}

Topological graph theory is useful for describing the topological structure of functions and flows on surfaces. These structures are given by discrete structures -- graphs with additional information. These graphs are often embedded in a surface, the resulting invariants are discrete Morse functions. There is a natural equivalence on the set of such functions, first explored in the papers of Forman \cite{for98}. At the same time, the equivalence class is given by a discrete gradient vector field. In the most general case, such fields and functions are considered on regular cellular complexes. The most useful discrete functions and vector fields are those that have the smallest number of critical cells in a given space, among all possible functions and vector fields on the regular cellular decomposition with minimal number of cells. Such fields are called optimal. Therefore, it is important to study the structure of such fields.

Graphs as topological invariants of functions were used in the papers of Kronrod \cite{Kronrod1950} and Reeb \cite{Reeb1946} for oriented maniofolds, in \cite{lychak2009morse} for  non-orientable two-dimensional manifolds and in   \cite{Bolsinov2004, hladysh2017topology, hladysh2019simple} for manifolds with boundary, in \cite{prishlyak2002morse} for non-compact manifolds. 

In general, Morse vector fields (Morse-Smale vector fields without closed orbits) are gradient field of Morse functions. If we fix the value of functions in singular points the field determinate the topological structure of the function \cite{lychak2009morse, Smale1961}. Therefore, Morse--Smale vector fields classification is closely related to the classification of the functions.

Topological classification of smooth function on closed 2-manifolds was also investigated in \cite{ hladysh2019simple, hladysh2017topology,  prishlyak2002morse, prishlyak2000conjugacy,  prishlyak2007classification, lychak2009morse, prishlyak2002ms, prish2015top, prish1998sopr,  bilun2002closed,  Sharko1993}, on 2-manifolds with the boundary in \cite{hladysh2016functions, hladysh2019simple} and on closed 3-manifolds in  \cite{prishlyak1999equivalence}.

In \cite{Kybalko2018, Oshemkov1998, Peixoto1973, prishlyak1997graphs, prishlyak2020three, akchurin2022three, prishlyak2022topological, prishlyak2017morse,  kkp2013,  prish2002vek,  prishlyak2021flows,  prishlyak2020topology,   prishlyak2019optimal}, 
the classifications of flows on closed 2- manifolds and 
\cite{loseva2016topology, prishlyak2017morse, prishlyak2022topological, prishlyak2003sum, prishlyak2003topological, prishlyak2019optimal} on manifolds with the boundary were obtained.
Topological properties of Morse-Smale vector fields on 3-manifolds was investigated in \cite{prish1998vek,  prish2001top, Prishlyak2002beh2, prishlyak2002ms,   prish2002vek, prishlyak2005complete, prishlyak2007complete, hatamian2020heegaard, bilun2022morse, bilun2022visualization}.
%prishlyak2003regular,

The main invariants of graphs and their embeddings in surfaces can be founded in \cite{prishlyak1997graphs, Harary69, HW68, GT87}.

%The topological properties of the projective plane are described in many geometry textbooks, but we would also recommend \cite{Bilun22Projective}.

The main purpose of this paper is to study the topological properties of the optimal discrete gradient vector fields on surface with one or two critical cells. These fields have and minimal number of critical cell on the regular cellular complex with minimal number of cells.

The first section gives the basic definitions and topological properties of discrete Morse functions and discrete gradient vector fields.

In the second section we describe all possible structures of optimal discrete vector flows on a two-dimensional disk, in the third section -- on a two-dimensional sphere, in the fourth section -- on the cylinder, and in the fifth section -- on a Mobius strip.

\section{Discrete Morse functions and vector field }

If a cell (or simplicial) complex $K$ is given, then a discrete function is a mapping that matches each simplex with a real number.
If the structure of a cellular complex (triangulation) is given on the manifold $M$, then for the function $f: M \to \mathbb{R}$ we take the value in the center of each cell (simplex). Thus, we  construct a discrete Morse function. On the contrary, if a Morse function is given on a triangulated manifold, then it determines the value of the function at the center of each simplex. Consider the first barycentric sub-division. For it, we have the value of the function at each vertex. By linearity, we extend the function to each simplex, thereby obtaining a piecewise linear function on the manifold. With the help of smoothing, you can get a smooth function from it. Similar constructions can be made for regular cage complexes.

Forman approach based on the use of discrete gradient flows was proposed in \cite{for98}.These designs are easier to build algorithms and implement on computers.

On the set of simplexes (cells) of the $K$ complex, we introduce the order relation:
$$\tau < \sigma \Leftrightarrow \tau \subset \sigma, \tau \ne \sigma.$$

By $|A|$ we denote the number of elements of the set $A$.

%\begin{definition} \index{discrete Morse function} 
The function $f: K \to \mathbb{R}$ is called a \textit{discrete Morse function} if for an arbitrary simplex (cell) $\sigma$ of dimension $k$:

$$ |\tau^{k-1}: \tau^{k-1}<\sigma, f(\tau^{k-1}) \ge f(\sigma)| +
|\tau^{k+1}: \tau^{k+1}>\sigma, f(\tau^{k+1}) \le f(\sigma)| \le 1.$$
%\end{definition}
%\begin{definition} \index{critical simplex} \index{critical cell} A simplex (cell) $\sigma$ of dimension $k$ is called critical for the function $f: K \to \mathbb{R}$ if:

$$ |\tau^{k-1}: \tau^{k-1}<\sigma, f(\tau^{k-1}) \ge f(\sigma)| +
|\tau^{k+1}: \tau^{k+1}>\sigma, f(\tau^{k+1}) \le f(\sigma)| = 0.$$

Simplexes that are not critical are called regular.
%\end{definition}

A discrete Morse function is called simple if it takes different values on different critical simplexes (cells).

%\begin{definition} \index{discrete vector field}
A discrete vector field on the complex $K$ is the set of pairs $V=\{ (\sigma, \tau) \}$ for which $\sigma, \tau \in K$, $\sigma <\tau$, $\ dim \sigma = \dim \tau +1$ and each simplex (cell) of the complex can belong to no more than one pair. Each such pair is called a vector and is represented by an arrow directed from the center of the first simplex (cell) to the center of the second.
%\end{definition}

%\begin{definition} \index{gradient field}

The gradient vector field $V_f$ of the discrete Morse function
is called the field $$V_f=\{ \sigma^{(k)}, \tau^{(k+1)}: \sigma^{(k)}<\tau^{(k+1)}, f( \sigma^{(k)}) \ge f( \tau^{(k+1)}).$$
%\end{definition}

\begin{figure}[ht!]
\begin{minipage}[h]{0.4\linewidth}
\center{\includegraphics[width=1\linewidth ]{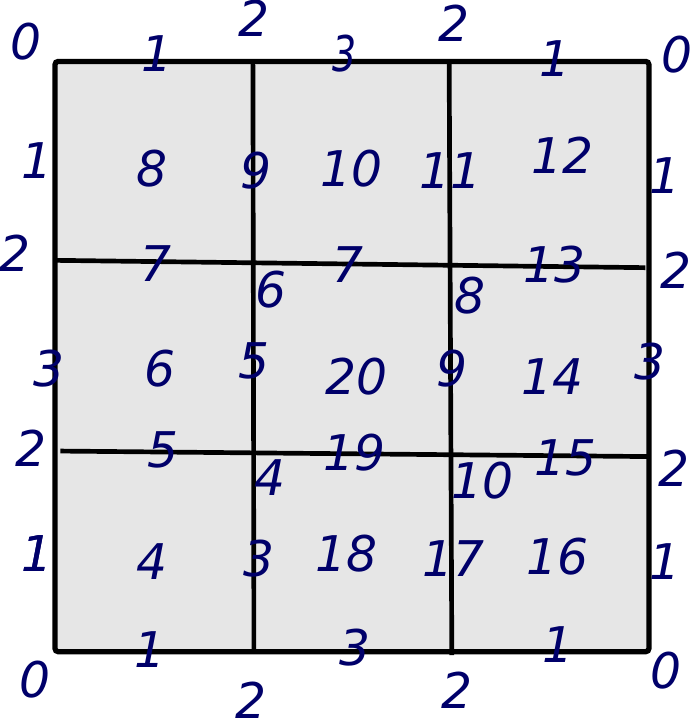}} a) \\
\end{minipage}
\hfill
\begin{minipage}[h]{0.4\linewidth}
\center{\includegraphics[width=1\linewidth]{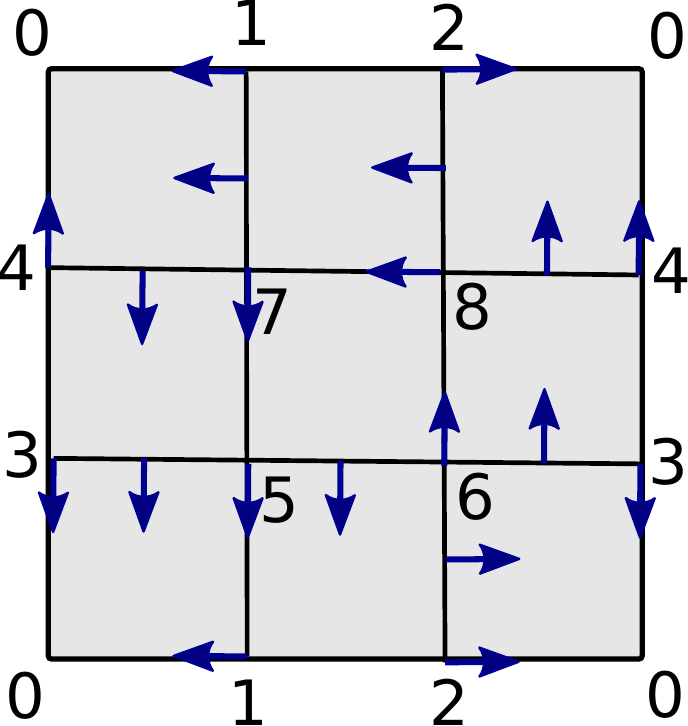}} \\b)
\end{minipage}
\caption{Discrete Morse function a) and its gradient field b) on the torus}
\label{dmft2}
\end{figure}

In fig. \ref{dmft2} shows the discrete Morse function on the torus and its gradient field. Critical cells are: 1) vertex 0, 2) edges $\{1,2\}$ and $\{3,4\}$, 3) 2- cell $[5,6,8,7]$.

The index of a critical cell (simplex) is its dimension.

%\begin{definition} \index{gradient path}
\textit{The gradient path} of a discrete vector field $V$ is called a sequence of cells (simplexes)
$$ \tau_0, \sigma_1, \tau_1, \ldots, \sigma_{m-1} , \tau_{m-1},\sigma_m$$
such that $$(\sigma_i, \tau_i)\in V, \sigma_{i}< \tau_{i-1}, 1 \le i \le m, \dim \sigma_0=\dim \sigma_m.$$
%\end{definition}

Examples of gradient paths on \ref{dmft2} b): 1) two paths from $\{1,2\}$ to $\{0\}$ -- $\{1,2\},\{1\} ,\{0,1\},\{0\}$ and $\{1,2\},\{2\},\{0,2\},\{0\}$; 2) two paths from $[5,6,8,7]$ to $\{1,2\}$ -- $[5,6,8,7],\{5,6\},[1, 2,6,5],\{1,2\}$ and $[5,6,8,7],$ $\{5,6\},$ $[1,2,6,5],\{2,6 \}$, \ $[2,0,3,6],\{3,6\},[6,3,4,8],\{8,4\},[8,4,0,2 ],\{2,8\},[8,2,1,7],\{2,1\}$.

%We introuduce next order on the set of all cell:
$$\tau < \sigma \Leftrightarrow \tau \subset \sigma, \tau \ne \sigma.$$

%A discrete vector field on the complex $K$ is the set of pairs $V=\{ (\sigma, \tau) \}$ for which $\sigma, \tau \in K$, $\sigma <\tau$, $\ dim \sigma = \dim \tau +1$ and each cell of the complex can belong to no more than one pair. Each such pair is called a vector and is represented by an arrow directed from the center of the first cell to the center of the second.

Two  discrete vector field are called isomorphic if there is an isomorphism of the CW-complexes which preserve the pairs.

Our aim  is to describe all possible up to isomorphism structures of flows with one or two critical cell. 

\section{Discrete vector fields on the 2-disk}
\begin{figure}[ht!]
\center{\includegraphics[width=0.15\linewidth ]{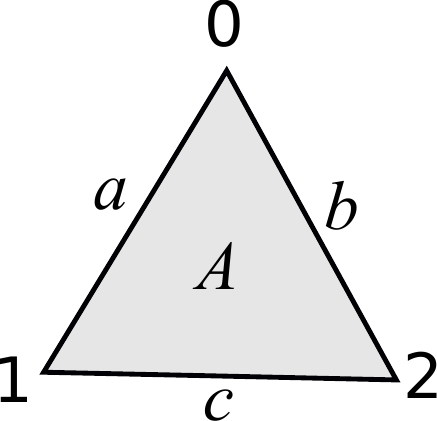}} 
\caption{the numaration of cells on the 2-disk}
\label{D2}
\end{figure}

Since a 2-disk is homotopically equivalent to a point, any discrete gradient vector field on it will have a critical 0-cell. The optimal cellular partition of the 2-disk is a triangle - 3 vertices, 3 edges and one face (it cannot be less, since there must be at least three vertices and three edges on the boundary of the surface to fulfill the definition of a regular cellular complex). To specify discrete vector fields on a 2-disk, we use the numbering and notation of cells, as in Fig. \ref{D2}.

Since the 2-cell is not critical, it belongs to one pair. Without loss of generality, we can assume that this is a pair of cA. Note that there is an axial symmetry that preserves this pair. Taking it into account, there are two possibilities for choosing a critical 0-cell: 0 or 1. In the first case, we will find such a vector field: 
$$V_1=\{1a, 2b, cA\}.$$ In the second case, the vector field is $$V_2=\{0a, 2b, cA\}.$$

Thus, the following is true.
\begin{theorem}
The number of non-isomorphic optimal discrete vector fields on the 2-disk is equal 2.
\end{theorem}

\section{Discrete vector fields on the 2-sphere}

\begin{figure}[ht!]
\center{\includegraphics[width=0.28\linewidth ]{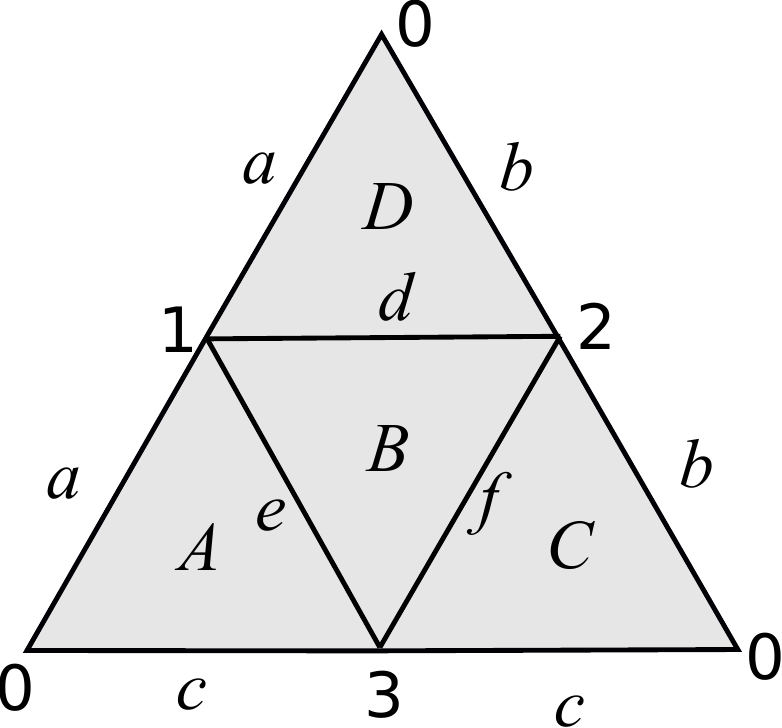}} 
\caption{the numaration of cells on the 2-shpere}
\label{S2}
\end{figure}

According to the discrete Morse theory, an optimal discrete gradient vector field on a 2-sphere has one critical 0-cell and one critical 2-cell. An optimal regular CW-complex of a sphere has four 0-cells, six 1-cells, and four 2-cells. Let us denote the cells of the regular cell partition as shown in Fig. \ref{S2}.

Without loss of generality, we assume that the critical 2-cell is cell $B$. If none of the 1-cells ($e, d, f$) bordering it are in pairs with 2-cells, then they form a cycle, which is impossible. Thus, three cases are possible: 1) one of these cells is paired with a 2-cell ($eA$), 2) two pairs ($eA$ and $fC$), 3) three pairs ($eA$, $fC$ and $dD$). In the first case, there are 2 variant:
$$V_{1-3}=\{\ldots, eA, aD, cC\},$$
$$V_{4-7}=\{\ldots, eA, aD, bC\}.$$

By $V_{1-3}$ we denote next three vector fields: 

In the first variant it is possible to select the critical vertex in 3 ways (due to symmetry):
$$V_{1}=\{2b,1d,3f, eA, aD, cC\},$$
$$V_{2}=\{0b,2d,3f, eA, aD, cC\},$$
$$V_{3}=\{0b,1d,3f, eA, aD, cC\}.$$ 
In the second variant we have 4 vector fields:
$$V_{4}=\{1d,2f,3c, eA, aD, bC\},$$
$$V_{5}=\{0c,2d,3f, eA, aD, bC\},$$
$$V_{6}=\{0c,1d,3f, eA, aD, bC\},$$
$$V_{7}=\{0c,1d,2f, eA, aD, bC\}.$$
In the second case $$V_{8-11}=\{\ldots, eA, fC, aD\}$$ there are four ways for critical vertex location. 

In the third case $$V_{12-13}=\{\ldots, eA, dD, fC\}$$ there are two ways for critical vertex location.
Thus, the following is true.
\begin{theorem}
The number of non-isomorphic optimal discrete vector fields on the 2-sphere is equal 13.
\end{theorem}

\section{Discrete vector fields on the cylinder } %$ S^1 \times [0,1] $ }
\begin{figure}[ht!]
\center{\includegraphics[width=0.3\linewidth ]{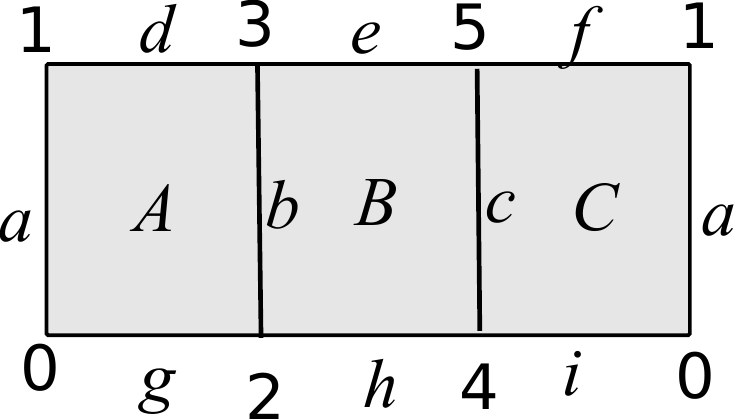}} 
\caption{the numaration of cells on the cylinder}
\label{SI}
\end{figure}

The optimal structure of the cell complex on a cylinder consists of three 2-cells, quadrilaterals, glued in a circle. All vertices in it are equal. Let's assign number 0 to the critical vertex. We cut the cylinder along the glue, ending at this vertex. we get a rectangle divided into three squares. We number its cells as in Fig. \ref{SI}.
Note that an optimal vector field has one critical vertex and a critical edge. There is one non-trivial isomorphism of the cell complex onto itself, which is fixed at the critical vertex. On fig. 3 it corresponds to axial symmetry about the vertical axis.

The vector field is determined by the pairs that contain the 2-cell and the location of the critical 1-cell. The following options are possible:
	 $$V_{1-3}=\{\ldots, aA, bB, iC\}, \
	 V_{4-6}=\{\ldots, aA, bB, fC\}, \
	 V_{7-9}=\{\ldots, gA, bB, aC\},  \
	 V_{10-12}=\{\ldots, gA, bB, iC\},$$
	$$V_{13-15}=\{\ldots, gA, bB, cC\}, \ \
	V_{16-20}=\{\ldots, gA, bB, fC\}, \ \
	V_{21-23}=\{\ldots, dA, bB, aC\},$$
	 $$V_{24-28}=\{\ldots, dA, bB, iC\}, \ \
	 V_{29-31}=\{\ldots, dA, bB, cC\}, \ \
	 V_{32-34}=\{\ldots, dA, bB, fC\},$$
	 $$V_{35-39}=\{\ldots, gA, eB, aC\}, \ \
	 V_{40-44}=\{\ldots, gA, eB, iC\}, \ \
	 V_{44-49}=\{\ldots, gA, eB, cC\},$$
	 $$V_{50-54}=\{\ldots, gA, eB, fC\}, \ \
	 V_{55-57}=\{\ldots, aA, eB, cC\},\ \
	 V_{58-60}=\{\ldots, aA, eB, fC\},$$
	 $$V_{61-63}=\{\ldots, bA, eB, cC\},\ \ 
	 V_{64-66}=\{\ldots, bA, eB, fC\},\ \
	 V_{67-69}=\{\ldots, dA, eB, fC\},$$
	 $$V_{70-72}=\{\ldots, gA, hB, aC\}, \ \
	 V_{73-75}=\{\ldots, gA, hB, iC\}, \ \
	 V_{76-78}=\{\ldots, gA, hB, cC\},$$
	 $$V_{79-83}=\{\ldots, gA, hB, fC\}, \ \ 
	 V_{84-86}=\{\ldots, aA, hB, cC\}, \ \
	 V_{87-91}=\{\ldots, aA, hB, fC\},$$
	 $$V_{92-94}=\{\ldots, bA, hB, cC\},\ \
	 V_{95-99}=\{\ldots, bA, hB, fC\}, \ \
	 V_{100-104}=\{\ldots, dA, hB, fC\}.$$

Thus, the following is true.
\begin{theorem}
The number of non-isomorphic optimal discrete vector fields on the cylinder $S^1\times [0,1]$  is equal 104.
\end{theorem}

\section{Discrete vector fields on the Mobius band}

The case with the Mobius band is similar to the case with the cylinder. The axial symmetry of the rectangle is replaced by central symmetry. Cell designations are shown in fig. 4.
\begin{figure}[ht!]
\center{\includegraphics[width=0.3\linewidth ]{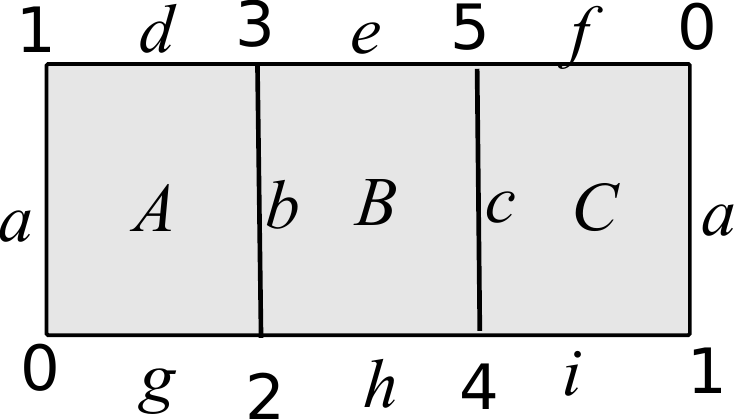}} 
\caption{the numaration of cells on the Mobius band}
\label{Mo}
\end{figure}

The vector field is determined by the pairs that contain the 2-cell and the location of the critical 1-cell. The following options are possible:
	 $$V_{1-4}=\{\ldots, aA, bB, iC\}, \
	 V_{5-8}=\{\ldots, aA, bB, fC\}, \
	 V_{9-12}=\{\ldots, gA, bB, aC\},  \
	 V_{13-16}=\{\ldots, gA, bB, iC\},$$
	$$V_{17-20}=\{\ldots, gA, bB, cC\}, \ \
	V_{21-24}=\{\ldots, gA, bB, fC\}, \ \
	V_{25-28}=\{\ldots, dA, bB, aC\},$$
	 $$V_{29-32}=\{\ldots, dA, bB, iC\}, \ \
	 V_{33-36}=\{\ldots, dA, bB, cC\}, \ \
	 V_{37-40}=\{\ldots, dA, bB, fC\},$$
	 $$V_{41-44}=\{\ldots, gA, eB, aC\}, \ \
	 V_{45-50}=\{\ldots, gA, eB, iC\}, \ \
	 V_{51-54}=\{\ldots, gA, eB, cC\},$$
	 $$V_{55-58}=\{\ldots, gA, eB, fC\}, \ \
	 V_{59-62}=\{\ldots, aA, eB, cC\},\ \
	 V_{63-66}=\{\ldots, aA, eB, iC\},$$
	 $$V_{58-70}=\{\ldots, aA, eB, fC\}, \ \
	 V_{71-74}=\{\ldots, bA, eB, cC\},\ \ 
	 V_{75-78}=\{\ldots, bA, eB, iC\},$$
	 $$V_{79-82}=\{\ldots, bA, eB, aC\},\ \
	 V_{83-86}=\{\ldots, bA, eB, fC\},\ \
	 V_{87-90}=\{\ldots, dA, eB, aC\}, $$
	$$V_{91-94}=\{\ldots, dA, eB, iC \}, \ \
	V_{95-98}=\{\ldots, dA, eB, fC\}, \ \
	V_{99-102}=\{\ldots, dA, eB, cC\}.$$

Thus, the following is true.
\begin{theorem}
The number of non-isomorphic optimal discrete vector fields on the cylinder $S^1\times [0,1]$  is equal 102.
\end{theorem}

\section*{Conclusion}
The results of the article prove the effectiveness of the proposed invariants (distinguishing graph and code) for classifying flows on a two-dimensional disk with one singular point. We hope that they can be used to construct similar invariants for flows with a large number of singular points on other surfaces as well.

%\bibliographystyle{plain}
%\bibliography{prishe}

\textsc{Taras Shevchenko National University of Kyiv}

Svitlana Bilun  \ \ \ \ \ \ \ \ \ \ \ \
\textit{Email address:} \text{ bilun@knu.ua}   \ \ \ \ \ \ \ \ \
\textit{ Orcid ID:}  \text{0000-0003-2925-5392}

Maria Hrechko \ \  \ \ \ \ \ \ \ \ \  \textit{Email:} \text{mariyagr2808@gmail.com }   \ \ \ \
\textit{ Orcid ID:} \text{0009-0005-5347-0342}

Olena Myshnova \ \ \ \ \ \ \ \ \
\textit{Email:} \text{ myshnova.olena@gmail.com}  \ 
\textit{ Orcid ID:} \text{0009-0001-8808-0038}

Alexandr Prishlyak \ \ \ \ \ 
\textit{Email address:} \text{ prishlyak@knu.ua} \ \ \ \ 
\textit{ Orcid ID:} \text{0000-0002-7164-807X}
\end{document}